\documentclass{amsart}

\usepackage{amsmath,amssymb,amscd,amsfonts}

\newtheorem{theorem}{Theorem}
\newcommand{\bt}{\begin{theorem}}
\newcommand{\et}{\end{theorem}}
\newtheorem{lemma}{Lemma}
\newcommand{\bl}{\begin{lemma}}
\newcommand{\el}{\end{lemma}}
\newtheorem{corollary}{Corollary}
\newcommand{\bc}{\begin{corollary}}
\newcommand{\ec}{\end{corollary}}
\newcommand{\beq}{\begin{equation}}
\newcommand{\eeq}{\end{equation}}
\newcommand{\benum}{\begin{enumerate}}
\newcommand{\eenum}{\end{enumerate}}
\newcommand{\N}{\ensuremath{ \mathbf N }}

\newcommand{\R}{\ensuremath{\mathbf R}}

\DeclareMathOperator{\Bad}{\text{Bad}}

\begin{document}

\title[Sets without geometric progressions]
{A problem of Rankin on sets without geometric progressions}
\author{Melvyn B. Nathanson}
\address{Department of Mathematics\\
Lehman College (CUNY)\\
Bronx, NY 10468}
\email{melvyn.nathanson@lehman.cuny.edu}

\author{Kevin O'Bryant}
\address{Department of Mathematics\\
College of Staten Island (CUNY)\\
Staten Island, NY 10314}
\email{kevin@member.ams.org}


\date{\today}

\begin{abstract}
A {\em geometric progression of length $k$ and integer ratio} is a set 
of numbers of the form $\{a,ar,\dots,ar^{k-1}\}$ 
for some positive real number $a$ and integer $r\geq 2$.   
For each integer $k \geq 3$, a greedy algorithm is used to construct 
a strictly decreasing  sequence 
$(a_i)_{i=1}^{\infty}$ of positive real numbers with $a_1 = 1$ such that the set 
\[
G^{(k)} = \bigcup_{i=1}^{\infty} \left(a_{2i} ,  a_{2i-1} \right]
\]
contains no geometric progression of length $k$ and integer ratio.
Moreover, $G^{(k)}$ is a maximal subset of $(0,1]$ that  
contains no geometric progression of length $k$ and integer ratio.
It is also proved that there is a strictly increasing  sequence 
$(A_i)_{i=1}^{\infty}$ of positive integers with $A_1 = 1$ 
such that $a_i = 1/A_i$ for all $i = 1,2,3,\ldots$.

The set $G^{(k)}$  gives a new lower bound for the maximum 
cardinality of a subset of 
the set of integers $\{1,2,\dots,n\}$ that contains no geometric 
progression of length $k$ and integer ratio.
\end{abstract}

\maketitle

\section{Real and integral geometric progressions}

Let \R\ denote the real numbers.  
For $t \in \R$, let $\R_{>t}$ denote the set of all real numbers $x > t$.  
Let $[x]$ denote the integer part of the real number $x$.  
For real numbers $u < v$,  we define the intervals 
\[
(u,v] = \{x\in \R: u < x \leq v \} \quad\text{and}\quad [u,v) = \{x\in \R: u \leq x < v \}.
\]

Let  $X$ be a set of positive real numbers, and let $u, v \in \R_{>0}$ with $u < v$.  
The \emph{dilation of the set $X$ by $q \in \R_{>0}$} 
is the set 
\[
q \ast X = \{qx : x \in X\}.  
\] 
The \emph{reciprocal} of the set $X$ is the set 
\[
X^{-1} = \left\{  x^{-1} : x \in X \right\}.
\]
For example, $q\ast (u,v] = (qu,qv]$ and 
$\left( 1/v, 1/u \right]^{-1} = [u,v)$.

If $A = (a_0,a_1,\ldots, a_{k-1})$ is a finite sequence of positive real numbers, 
then the \emph{dilation of the sequence $A$ by $q$}  
is the sequence $q \ast A = (qa_0,qa_1,\ldots, qa_{k-1})$ 
and the reciprocal of $A$ is the sequence  
$A^{-1} = (1/a_0, 1/a_1,\ldots, 1/a_{k-1})$.

Let \N\ denote the set of positive integers, 
and let $\N^{\sharp} = \N\setminus \{1\}$ 
denote the set of all integers $r > 1$.  
Let $k \in \N$ and let $r,a \in \R_{>0}$.  
A \emph{geometric progression of length $k$ and ratio $r$ with first term $a$} 
is a sequence of the form 
\[
 (a, ar, ar^2, \ldots, ar^{k-1} ) = a \ast (1,r,r^2,\ldots, r^{k-1}).
\]
This is an \emph{integer geometric progression of length $k$}
if $a r^j  \in \N$ for all $j \in \{0,1,\ldots, k-1\}$.
If $(a, ar, ar^2, \ldots, ar^{k-1} )$ is an integer geometric progression, then 
the ratio $r$ must be a rational number.  
For example, $(8, 12,18, 27)$ is an integer geometric progression of length 4 
with ratio 3/2.

Note that the dilation  by a positive real number $q$ 
of the geometric progression  $(a, ar, ar^2, \ldots, ar^{k-1} )$ 
of length $k$, ratio $r$, and first term $a$ is the 
geometric progression   $(qa, qar, qar^2, \ldots, qar^{k-1} )$  
of length $k$, ratio $r$, and first term $qa$.   
The reciprocal of the geometric progression  $(a, ar, ar^2, \ldots, ar^{k-1} )$ 
is the geometric progression  
\[
\left( \frac{1}{a}, \frac{1}{a}\left( \frac{1}{r} \right), \frac{1}{a}\left( \frac{1}{r} \right)^2, 
\ldots, \frac{1}{a}\left( \frac{1}{r} \right)^{k-1} \right)
\] 
of length $k$, ratio $1/r$, and first term $1/a$.   

The \emph{reverse} of  the sequence $(a_1, a_2,\ldots, a_{k-1}, a_k)$ 
is the sequence $(a_k, a_{k-1}, \ldots, a_2,a_1)$.  
The reverse of the reciprocal of the geometric progression  $(a, ar, ar^2, \ldots, ar^{k-1} )$ is the geometric progression  $(b, br, br^2, \ldots, br^{k-1} )$, 
where $b = 1/\left(ar ^{k-1}\right)$.

Thus, a set $G$ of real numbers  contains no geometric progression of length $k$ 
if and only if the dilation $q\ast G$ contains no geometric progression of length $k$ 
for every positive real number $q$.  
Moreover, if a set $G$ contains no geometric progression of length $k$, 
then no subset of $G$ contains a geometric progression of length $k$.
It follows that if a set $G$ contains no geometric progression of length $k$, 
then, for every positive real number $q$,  the set $(q \ast G)\cap \N$ 
is a set of positive integers that contains no geometric progression of length $k$.  
Similarly, if $G$ contains no geometric progression of length $k$, 
then the set of $G^{-1}\cap \N$ 
is a set of positive integers that contains no geometric progression of length $k$.  

A \emph{geometric progression of length $k$  with integer ratio} 
is a geometric progression of length $k$  with ratio $r \in \N^{\sharp}$.
An \emph{integer geometric progression of length $k$  with integer ratio} 
is a geometric progression of the form $(a,ar, ar^2, \ldots, ar^{k-1} )$  
with $a \in \N$ and $r \in \N^{\sharp}$.

For positive integers $k$ and $n$, let $g_k(n)$ denote 
the cardinality of the largest subset of the set $\{1,2,3,\ldots, n\}$ 
that contains no integer geometric progression of length $k$ with integer ratio,  
and  let $\hat{g}_k(n)$ denote 
the cardinality of the largest subset of the set $\{1,2,3,\ldots, n\}$ 
that contains no integer geometric progression of length $k$ with rational ratio.

We have $g_1(n) = \hat{g}_1(n) = 0$ for all $n \in \N$, 
and $g_k(n) =  \hat{g}_k(n) = n$ if $n < k$.  
Moreover, $\hat{g}_2(n) = 1$ for $n \geq 2$.  
We compute $g_2(n)$ in the next section.  
In this paper we obtain new lower bounds for the function $g_k(n)$ for $k \geq 3$.    

For every integer $k \geq 3$, there are four basic unsolved problems:
\benum
\item
Determine the cardinality and the structure of the maximal subsets of $\{1,2,\ldots, n\}$ 
that contain no geometric progression of length $k$ with integer ratio.  
In particular, what is the maximum cardinality $g_k(n)$?

\item
Determine the cardinality and the structure of the maximal subsets of $\{1,2,\ldots, n\}$ 
that contain no geometric progression of length $k$ with rational ratio.
What is the maximum cardinality $\hat{g}_k(n)$?

\item
Determine the density and  structure of maximal infinite sets of positive integers  
that contain no geometric progression of length $k$ with integer ratio.
What is the least upper bound of the densities of such sets? 
Is this least upper bound achieved?

\item
Determine the density and  structure of maximal infinite sets of positive integers  
that contain no geometric progression of length $k$ with rational ratio.
What is the least upper bound of the densities of such sets? 
Is this least upper bound achieved?
\eenum

Very little is known about these problems.  The literature consists mostly of lower bounds 
for the maximum cardinalities in Problems 1 and 2, 
and for the densities in Problems 3 and 4.  
In this paper we improve the  lower bounds in Problem 1.  
Our method is to use a greedy algorithm to construct, 
for every integer $k \geq 3$,  a unique maximal subset of 
the unit interval $(0,1]$ that contains no geometric progression of length $k$ 
with integer ratio,
and to use the measure of this set to obtain new lower bounds for the finite sets 
considered in Problem 1.  

The earliest discussion of sets with no $k$-term geometric progression 
is in a paper of Rankin~\cite{rank60} in 1960 that 
was concerned with sets of integers containing no $k$-term arithmetic progression.

\section{Integral geometric progressions of length 2}
We can quickly solve the problem of integer geometric progressions 
of length 2 with integer ratio.  
Every set $\{a,b\}$  of positive real numbers with $a < b$ 
is a geometric progression of length 2 with ratio $r = b/a$.  
In particular, every set $\{a,b\}$  of positive integers with $a < b$ 
is a geometric progression of length 2 with rational ratio $r = b/a$.  
The set $\{a,b\}$ is an integer geometric progression of length 2 with integer ratio 
if and only if $a,b \in \N$ and $a$ divides $b$.  
Thus, a set $S$ of positive integers contains no 2-term geometric progression 
if and only if $S$ is \emph{primitive} in the sense that no element of $S$ 
divides another element of $S$.  

The following is a classical result in combinatorial number theory.

\bt
Let $g_2(n)$ denote the cardinality of the largest primitive 
subset of $\{1,2,\ldots, n\}$, that is, 
the largest subset of $\{1,2,\ldots, n\}$ that contains 
no  integer geometric progression of length 2 with integer ratio. 
Then $g_2(n) =  \left[\frac{n+1}{2}\right]$. 
\et

\begin{proof}
For every positive integer $n$, the interval 
\beq                      \label{NoGPunit:primitive} 
S = \left(  \left[\frac{n}{2}\right], n \right] 
= \left\{  \left[\frac{n}{2}\right]+1, \left[\frac{n}{2}\right] +2, \ldots, n -1, n \right\}
\eeq
is primitive because $2\left( \left[\frac{n}{2}\right]+1 \right) \geq n+1 > n$.  
The cardinality of this set is $ \left[\frac{n+1}{2}\right]$, 
and so $g_2(n) \geq \left[\frac{n+1}{2}\right]$.

Let $S$ be any primitive subset of $\{1,2,\ldots, n\}$.  
Each element $s \in S$ can be written uniquely in the form $s = 2^{k(s)} a(s)$, 
where $k(s)$ is a nonnegative integer and $a(s)$ is an odd integer in $\{1,2,\ldots, n\}$.
If $a(s_1) = a(s_2)$ for integers $s_1, s_2 \in S$ with $s_1 < s_2$,   
then $s_1$ divides $s_2$.  It follows that the cardinality of the primitive set $S$ 
is at most the number of odd integers in $\{1,2,\ldots, n\}$, 
and so $g_2(n) \leq \left[\frac{n+1}{2}\right]$.
This completes the proof.  
\end{proof}

\section{Good set, bad set}
Let $k$ be an integer, $k \geq 3$.   
A \emph{$k$-good set} is a set of positive real numbers that contains 
no geometric progression of length $k$ with integer ratio.  
For example, the set
\[
G_1^{(k)} = \left( \frac{1}{2^{k-1}}, 1 \right]
\]
is $k$-good because, if $x\in G_1^{(k)}$ and $r\in \N^{\sharp}$, 
then $xr^{k-1} \geq x2^{k-1} > 1$ 
and so $\{x, xr, xr^2, \ldots, xr^{k-1}  \}$ is not a subset of $G_1^{(k)}$.

Let $G$ be a $k$-good subset of $(0,1]$, and let $x \in (0,1] \setminus G$.  
The real number $x$ is \emph{$k$-bad with respect to  $G$} 
if there exists an integer $r \in \N^{\sharp}$ such that  
$G \cup \{ x \}$ contains the $k$-term geometric progression 
$(x, xr, xr^2, \ldots, xr^{k-1})$.  
Thus, if $x$ is $k$-bad with respect to $G$, 
then the set $G \cup \{x\}$ is not $k$-good.  

For example, the number $1/2^k$ is $k$-bad with respect to the 
$k$-good  set $G_1^{(k)}$ because 
$(1/2^k,1/2^{k-1}, 1/2^{k-2}, \ldots, 1/2, 1)$ 
is a $k$-term geometric progression with ratio $r=2$ 
contained in $G_1^{(k)} \cup \{1/2^k \}$.  
 
The number $3/16$ is $3$-bad with respect to the 3-good set 
$G_1^{(3)} = (1/4,1]$ because, with $r=2$,  
\[
\left\{ \frac{3}{16} r, \frac{3}{16} r^2\right\} = \left\{ \frac{3}{8}, \frac{3}{4} \right\} 
\subseteq \left( \frac{1}{4}, 1 \right] = G_1^{(3)}
\]
and so the set $G_1^{(3)}\cup \{ 3/16\}$ contains the 3-term geometric 
progression $(3/16, 3/8, 3/4)$.  
Similarly, $1/10$ is $3$-bad with respect to $G_1^{(3)}$ because, with $r=3$,  
\[
\left\{ \frac{1}{10} , \frac{1}{10} r, \frac{1}{10} r^2\right\} = \left\{\frac{1}{10},  \frac{3}{10}, \frac{9}{10} \right\} 
\subseteq G_1^{(3)} \cup \left\{ \frac{1}{10} \right\}.
\]

Note that if $G$ is a $k$-good subset of $(0,1]$ 
and if $x \in (0,1] \setminus G$ is $k$-bad with respect to $G$, 
then $x$ is also $k$-bad with respect to the good set $G \cap (x,1]$, 
because $x < r^jx$ for all $r \in \N^{\sharp}$ and $j \in \{1,\ldots, k-1\}$.

The real number $x  \in (0,1] \setminus G$ is \emph{$k$-good with respect to $G$} 
if $x$ is not $k$-bad with respect to $G$.  
Thus, $x$ is $k$-good with respect to $G$ if and only if, 
for every $r \in \N^{\sharp}$, there exists $j \in \{1,2,\ldots, k-1\}$ 
such that $ x r^j \notin G$.  
Because $G \subseteq (0,1]$ and $x r^{k-1} \notin G$ if $r > (1/x)^{1/(k-1)}$, 
it follows that  $x \in [0,1)\setminus G$ is $k$-good with respect to $G$ 
if and only if, for every integer $r$ with $2 \leq r \leq (1/x)^{1/(k-1)}$, 
there exists $j \in \{1,2,\ldots, k-1\}$ with $x r^j \notin G$.

For every $k$-good set $G \subseteq (0,1]$, we define 
\[
\Bad(G) = \{x\in (0,1] \setminus G :\text{ $x$ is $k$-bad with respect to $G$} \}.
\]
Thus, $G\cup \{x\}$ is $k$-good for all 
$x \in (0,1] \setminus \left( G \cup \Bad(G) \right)$.
If $G$ and $G'$ are $k$-good sets with $G\subseteq G'$, then 
$\Bad(G) \subseteq \Bad(G')$.

For fixed $k$, we usually write ``good'' instead of ``$k$-good'' 
and ``bad'' instead of ``$k$-bad.''

\section{Construction of a good set of real numbers}
Fix the integer $k \geq 3$.  
We shall use a greedy algorithm to construct a large good set contained 
in the  interval $(0,1]$.  
We begin with some simple observations about good and bad sets.  

\bl                    \label{NoGPunit:lemma:intervals} 
Let $k \geq 3$, let $0 < a < 1$, and let $\delta_k(a) = a^{(k-1)/(k-2)}$.  
\benum
\item[(i)]
For every $\delta >0$, the interval  $(0, \delta]$ is not good.
\item[(ii)]
Every number in the interval $(0,a^2]$ is good with respect to the interval $(a,1]$.
\item[(iii)]   
Let $x \in (0,1]$.  
If $xr^j \in (a,1]$ for some $r \in \N^{\sharp}$ 
and all $j \in \{1,\ldots, k-1\}$, then 
$x > \delta_k(a)$.

\item[(iv)]   If $G$ is a good set with $G \subseteq (a, 1]$, then $(0,\delta_k(a)] \cap \Bad(G) = \emptyset$.  
\eenum
\el

Note that $0 < \delta_k(a) < a$.

\begin{proof}
 (i) We have $0 < 2^{1-k} \delta < \delta$.  
For every $x \in (0, 2^{1-k} \delta]$, we have 
 \[
0 < x < 2x < \cdots <  2^{k-1}x  \leq 2^{k-1} 2^{1-k} \delta  = \delta
 \]
 and so $(0, \delta]$ contains the $k$-term geometric progression 
 $\{x,2x,\ldots, 2^{k-1}x \}$. 
Thus, the interval $(0,\delta]$ is not good.
 
 (ii)  
If $x \in (0, a^2]$, $r \in \N^{\sharp}$ and $xr\in (a,1]$, then $xr > a$ and $r > a/x$.
It follows that 
\[
xr^{k-1} \geq xr^2 > x\left(\frac{a}{x}\right)^2 = \frac{a^2}{x} > 1
\]
and so $x$ is good with respect to $(a,1]$.

$x$ is bad with respect to $(a,1]$, 
then there exists $r \in \N^{\sharp}$ such that $xr^i \in (a,1]$

 (iii)
 If $r \in \N^{\sharp}$ and 
 \[
a < rx < \cdots < r^{k-1}x \leq 1
 \]
then 
\[
\frac{a}{r} < x \leq \frac{1}{r^{k-1}}
\]
and so $1/r > a^{1/(k-2)}$.  Therefore, 
\[
x > \frac{a}{r} > a a^{1/(k-2)} = a^{(k-1)/(k-2)} = \delta_k(a).
\] 

(iv)  This follows immediately from~(iii).
\end{proof}

\bl                    \label{NoGPunit:lemma:OpenBad} 
Let  $(a_i)_{i=1}^{2n}$ be a strictly decreasing  sequence 
of positive real numbers with $a_1 \leq 1$ such that 
\[
G_n = \bigcup_{i=1}^n (a_{2i}, a_{2i-1}]
\]
is a good set.  
If $x \in \Bad(G_n)$, then there exists $\delta > 0$ such that 
$(x-\delta, x] \subseteq \Bad(G_n)$.
\el

\begin{proof}
We have $G_n \subseteq (0,1]$.  
If $x \in \Bad(G_n)$, then there exists $r \in \N^{\sharp}$ such that 
$x r^j \in G_n$ for all $j \in \{1,\ldots, k-1\}$.
It follows that, for each $j \in \{1,\ldots, k-1\}$, there exists $i_j \in \{1,\ldots, n\}$ 
such that 
\[
x r^j \in (a_{2i_j}, a_{2i_j -1}]
\]
or, equivalently,
\[
\frac{a_{ 2i_j}}{r^j} < x \leq \frac{ a_{2i_j-1} }{r^j}.
\]
Choose $\delta > 0$ such that 
\[
\frac{a_{ 2i_j}}{r^j} < x - \delta < x \leq \frac{ a_{2i_j-1} }{r^j}
\]
for all $j \in \{1,\ldots, k-1\}$.
If $y \in (x-\delta, x]$, then 
\[
a_{2i_j} <  (x-\delta)  r^j <   y r^j  \leq   xr^j \leq a_{2i_j-1} 
\]
and so $y r^j \in (a_{2i_j}, a_{2i_j-1}  ] \subseteq G_n$  
for all $j \in \{1,\ldots, k-1\}$.
Thus, $(x-\delta, x] \subseteq \Bad(G_n)$.
\end{proof}

\bl                   \label{NoGPunit:lemma:OpenGood}
Let  $(a_i)_{i=1}^{2n+1}$ be a strictly decreasing  sequence 
of positive real numbers with $a_1 \leq 1$ 
such that 
\[
G_n = \bigcup_{i=1}^n (a_{2i}, a_{2i-1}]
\]
is a good set, and   
\[
\bigcup_{i=1}^n (a_{2i+1}, a_{2i}] \subseteq \Bad(G_n).
\]
If $ x\in (a_{2n+1}/2, a_{2n+1}]$ is good with respect to $G_n$, 
then there exists $\delta > 0$ such that 
$(x - \delta, x] \cup G_n$ is good.  
\el

\begin{proof}
Let $x \in (a_{2n+1}/2, a_{2n+1}]$  be good with respect to $G_n$.  
For each $r \in \N^{\sharp}$ there exists $j_r \in \{1,\ldots, k-1\}$ 
such that $x r^{j_r} \notin G_n$.  
Let $r_0$ be the smallest integer such that $r_0 \geq 2$ and $x r_0^{k-1} > a_1$.  
Then $x > a_1/r_0^{k-1}$,  
and there exists $\delta_0 > 0$ such that $x - \delta_0 > a_1/r_0^{k-1}$.
If $y \in (x-\delta_0, x]$ and $r \geq r_0$, then 
\[
y r^{k-1}  >  (x-\delta_0) r_0^{k-1} > a_1
\]
and so $y r^{k-1} \notin G_n$.  

For each integer $r$ such that $2 \leq r < r_0$, 
we have 
\[
a_{2n+1} <  2x \leq x r^{j_r}  \leq x r^{k-1} \leq a_1
\]
and so there exists $i_r \in \{1,\ldots, n\}$ such that 
\[
a_{2i_r+1} < x r^{j_r} \leq a_{2i_r}.
\]
Equivalently, 
\[
\frac{a_{2i_r+1}}{r^{j_r}} < x \leq \frac{a_{2i_r}}{r^{j_r}}.
\]
Choose $0 < \delta_1 < x/2$ such that  
\[
\frac{a_{2i_r+1}}{r^{j_r}} < x - \delta_1 < x \leq \frac{a_{2i_r}}{r^{j_r}}  
\]
for all $r \in \N^{\sharp}$ with $r < r_0$.  
If $y \in (x-\delta_1, x]$ and $r < r_0$, then
\[
a_{2i_r+1} < (x-\delta_1) r^{j_r} <  y r^{j_r}\leq x r^{j_r}\leq a_{2i_r} 
\]
and so $y r^{j_r}  \notin G_n$.  
Let $\delta = \min(\delta_0,\delta_1)$.  
It follows that if $y \in (x-\delta, x]$, then  $y$ is good with respect to $G_n$.
This completes the proof.  
\end{proof}

\bt                   \label{NoGPunit:theorem:construct} 
Let $k \geq 3$.  
There exists  a unique strictly decreasing  sequence $(a_i)_{i=1}^{\infty}$ 
of positive real numbers with $a_1 = 1$ such that 
\[
G = \bigcup_{i=1}^{\infty} (a_{2i}, a_{2i-1}]
\]
is a good set, and 
\[
\Bad(G) =  \bigcup_{i=1}^{\infty} (a_{2i+1}, a_{2i}].
\]
\et

\begin{proof}
We construct the sequence $(a_i)_{i=1}^{\infty}$ by induction.

Let $a_1 = 1$.  If $x > 2^{1-k}$, then for all $r \in \N^{\sharp}$ we have
\[
r^{k-1}x > 2^{k-1} 2^{1-k} = 1
\]
and so $(2^{1-k},1]$ is a good set.  
Therefore, 
\[
a_2 = \inf\{ x\in (0,1] : \text{$(x,a_1] $ is good} \} \leq 2^{1-k}.
\]
We observe that $[2^{1-k},1]$ is not a good set because, 
with $y = 2^{1-k}$, we have  
$\{ y, y2,\ldots, y2^{k-1} \} \subseteq [2^{1-k},1]$.  Therefore, 
\[
a_2  = \frac{1}{2^{k-1}}\in \Bad(G_1)
\]
where
\[
G_1 = (a_2,a_1] = \left( \frac{1}{2^{k-1}},1 \right].
\]
We define 
\[
a_3 = \inf\{ x\in (0,1] : (x,a_2] \subseteq \Bad(G_1) \}.
\]
It follows from Lemmas~\ref{NoGPunit:lemma:intervals} 
and~\ref{NoGPunit:lemma:OpenBad} 
that $0 < \delta_k(a_2) \leq a_3  < a_2$ and $a_3 \notin \Bad(G_1)$.

Let $n \geq 1$, and assume that there is 
a unique strictly decreasing  sequence $(a_i)_{i=1}^{2n+1}$ 
of positive real numbers with $a_1 = 1$ such that 
\[
G_n = \bigcup_{i=1}^{n} (a_{2i}, a_{2i-1}]
\]
is a good set,
\[
\bigcup_{i=1}^{n} (a_{2i+1}, a_{2i}] \subseteq \Bad(G_n).  
\]
and
\[
a_{2n+1} = \inf\{ x\in (0,1] : (x,a_{2n}] \subseteq \Bad(G_n) \}.
\]
By Lemma~\ref{NoGPunit:lemma:OpenBad}, the number $a_{2n+1}$ 
is good with respect to $G_n$.  
Let 
\[
a_{2n+2} = \inf\{ x\in (0,a_{2n+1} ] : \text{$(x, a_{2n+1}]$ is good with respect to $G_n$}\}.
\]
Let 
\[
G_{n+1} = G_n \cup (a_{2n+2}, a_{2n+1}].
\]
Lemmas~\ref{NoGPunit:lemma:intervals} 
and~\ref{NoGPunit:lemma:OpenGood} imply that $0 < a_{2n+2} < a_{2n+1}$, 
and that $a_{2n+2} \in \Bad(G_{n+1})$.  
We define 
\[
a_{2n+3} = \inf\{ x\in (0,1] : (x,a_{2n+2}] \subseteq \Bad(G_n) \}.
\]
This completes the induction. 
\end{proof}

\bt                    \label{NoGPunit:theorem:harmonic-finite}
Let  $(a_i)_{i=1}^{2n}$ be a strictly decreasing  sequence 
of positive real numbers such that 
\[
G_n = \bigcup_{i=1}^n (a_{2i}, a_{2i-1}]
\]
is a good set, and 
\[
 \bigcup_{i=1}^{n-1} (a_{2i+1}, a_{2i}] \subseteq \Bad(G_n).
\]
If $A_1$ and $A_2$ are positive integers such that $a_1 = 1/A_1$ and 
$a_2 = 1/A_2$, then there is a strictly increasing  sequence  $(A_i)_{i =1}^{2n}$ 
of positive integers such that 
\[
a_i = \frac{1}{A_i}
\]
for $i = 1, \ldots,  2n$.
\et

\begin{proof}
The proof is by induction on $i$.    Let $2 \leq i \leq n$ and 
assume that there are positive integers $A_1 < \cdots < A_{2i-2}$ 
such that $a_j = 1/A_j$ for $j = 1,\ldots, 2i-2$.  
We shall prove that there are positive integers $A_{2i-1}$ 
and $A_{2i}$ such that $a_{2i-1} = 1/A_{2i-1}$ and  $a_{2i} = 1/A_{2i}$.

Consider the good number $a_{2i-1}$.  
If $h\in \N$  and $h \geq (a_{2i-2}-a_{2i-1})^{-1}$, then 
\[
a_{2i-1} + \frac{1}{h} \in (a_{2i-1},a_{2i-2}] \subseteq \Bad(G_n)
\]
and so there exists $r_h \in \N^{\sharp}$ such that, for all $j \in \{1,2,\ldots, k-1\}$, 
\[
\left( a_{2i-1} + \frac{1}{h}  \right) r_h^j  \in G_n 
\]
and
\[
a_{2i-1}   < a_{2i-1}  r_h  \leq   a_{2i-1} r_h^j 
< \left( a_{2i-1} + \frac{1}{h}  \right)  r_h^{k-1} \leq a_1. 
\]
Therefore, 
\[
2 \leq r_h < \frac{a_1}{a_{2i-1}}.
\]
Because $a_{2i-1} \in G_n$, there exists $j_h \in \{1,2,\ldots, k-1\}$ such that 
\[
 a_{2i-1} r_h^{j_h} \notin G_n. 
\]
There are only finitely many choices for $r_h$ and $j_h$.  
By the pigeonhole principle,  there are integers $r \in \N^{\sharp}$ 
and $j  \in \{1,2,\ldots, k-1\}$ 
and there is a strictly increasing infinite 
sequence $(h_{\ell})_{\ell \in \N}$ of positive integers 
such that 
\[
r_{h_{\ell}} = r \qquad\text{and} \qquad j_{h_{\ell}} = j
\]
for all $\ell \in \N$.  
Because $a_{2i-1} r^j  \notin G_n$ and $a_{2i-1} < a_{2i-1} r^j < a_1$, 
there is a unique positive integer $t \leq i$ 
such that $ a_{2i-1} r^j \in (a_{2t-1}, a_{2t-2} ]$.  
Because $ \left( a_{2i-1} + 1/h_{\ell}  \right) r^j  \in G_n$, 
it follows that  
\[
a_{2i-1}  r^j   \leq a_{2t-2} < \left( a_{2i-1} + \frac{1}{h_{\ell}}  \right)  r^j 
\]
or, equivalently,
\[
\frac{a_{2t-2}}{r^j } - \frac{1}{h_{\ell}}   < a_{2i-1}  \leq \frac{a_{2t-2}}{r^j }.
\]
By the induction hypothesis, there is a positive integer $A_{2h-2}$ such that 
$a_{2t-2} = 1/A_{2t-2}$.  
Letting $\ell \rightarrow \infty$, we obtain
\[
a_{2i-1}  = \frac{a_{2t-2}}{r^j } = \frac{1}{r^j A_{2t-2}} = \frac{1}{A_{2i-1}}
\]
with $A_{2i-1} = r^jA_{2t-2}$.

Next we consider the bad number $a_{2i}$.  
There exists $r \in \N^{\sharp}$ such that  $a_{2i}  r^j \in G_n$ 
for all $j \in \{1,2,\ldots, k-1\}$.  
If  $h \geq (a_{2i-1}-a_{2i})^{-1}$, then 
\[
a_{2i} + \frac{1}{h} \in (a_{2i},a_{2i-1}] \subseteq G_n
\]
and so there exists $j_h \in \{1,2,\ldots, k-1\}$ such that 
\[
\left( a_{2i} + \frac{1}{h}  \right) r^{j_h}  \notin G_n.
\]
By the pigeonhole principle,  there is an integer $j  \in \{1,2,\ldots, k-1\}$ 
and there is a strictly increasing infinite 
sequence $(h_{\ell})_{\ell \in \N}$ of positive integers 
such that $j_{h_{\ell}} = j$ for all $\ell \in \N$.  
Because $a_{2i} r^j \in G_n$, there is a unique positive integer $t \leq i$ such that 
$ a_{2i} r^j  \in (a_{2t}, a_{2t-1}]$.  
Because $\left( a_{2i} + 1/h_{\ell} \right) r^j  \notin G_n$, 
it follows that, for all $\ell \in \N$, we have 
\[
a_{2i}  r^j  \leq a_{2t-1} <  \left( a_{2i} + \frac{1}{h_{\ell}}  \right) r^j 
\]
or, equivalently,
\[
\frac{a_{2t-1}}{r^j } - \frac{1}{h_{\ell}}   < a_{2i}  \leq \frac{a_{2t-1}}{r^j }.
\]
By the induction hypothesis, there is a positive integer $A_{2t-1}$ such that 
$a_{2t-1} = 1/A_{2t-1}$. 
Letting $\ell \rightarrow \infty$, we obtain
\[
a_{2i}  = \frac{a_{2t-1}}{r^j } = \frac{1}{r^j A_{2t-1}} = \frac{1}{A_{2i}} 
\]
with $A_{2i} = r^j A_{2t-1}.$
This completes the proof.
\end{proof}

\bt          \label{NoGPunit:theorem:harmonic-infinite}
Let  $(a_i)_{i \in \N}$ be a strictly decreasing infinite sequence 
of positive real numbers such that 
\[
G = \bigcup_{i=1}^{\infty} (a_{2i}, a_{2i-1}]
\]
is a good set, and 
\[
\Bad(G) =  \bigcup_{i=1}^{\infty} (a_{2i+1}, a_{2i}].
\]
If $A_1$ and $A_2$ are positive integers such that $a_1 = 1/A_1$ and 
$a_2 = 1/A_2$, then there is a strictly increasing infinite sequence  $(A_i)_{i \in \N}$ 
of positive integers such that 
\[
a_i = \frac{1}{A_i}
\]
for all $i \in \N$.
Moreover, 
\[
\lim_{i \rightarrow \infty} a_i = 0.
\]
\et

\begin{proof}
Apply Theorem~\ref{NoGPunit:theorem:harmonic-finite} 
to the good set $G_n = \bigcup_{i=1}^n (a_{2i}, a_{2i-1}]$.

Because there is a strict increasing sequence $(A_i)_{i=1}^{\infty}$ of positive integers 
such that $a_i = 1/A_i$, it follows that 
\[
\lim_{i \rightarrow \infty} a_i = \lim_{i \rightarrow \infty} \frac{1}{A_i} = 0.
\]
This completes the proof.  
\end{proof}

\bt               \label{NoGPunit:theorem:constructHarmonic}
Let $k \geq 3$.  
There exists  a unique strictly increasing  sequence 
$\left(A_i^{(k)} \right)_{i=1}^{\infty}$ 
of positive integers with $A^{(k)}_1 = 1$ such that 
\beq            \label{NoGPunit:constructHarmonic}
G^{(k)} = \bigcup_{i=1}^{\infty} \left( \frac{1}{ A^{(k)}_{2i}} , \frac{1}{ A_{2i-1}^{(k)} } \right]
\eeq
is a $k$-good set and 
\[
\Bad\left( G^{(k)} \right) =  \bigcup_{i=1}^{\infty}  \left( \frac{1}{ A^{(k)} _{2i+1}} , \frac{1}{ A^{(k)} _{2i} } \right].
\]
\et

\begin{proof}
The existence and uniqueness of the sequence 
$\left(A_i^{(k)} \right)_{i=1}^{\infty}$  follows immediately from Theorems~\ref{NoGPunit:theorem:construct} 
and~\ref{NoGPunit:theorem:harmonic-infinite}.  
\end{proof}

Note that
\[
\inf G^{(k)}  = \inf \Bad\left( G^{(k)} \right)  = 0
\]
because  $\lim_{i\rightarrow \infty} A_i^{(k)} = \infty$.

We have already proved that $A_2^{(k)} = 2^{k-1}$.  
We can also  determine the integers $A_3^{(k)}$ and $A_4^{(k)}$.   
The proofs use a simple arithmetic inequality:
If $k \geq 3$, then 
\[
\frac{3^{k-1}}{2^k} = \frac{1}{2} \left( \frac{3}{2}\right)^{k-1} 
\geq \frac{1}{2} \left( \frac{3}{2}\right)^2 = \frac{9}{8} > 1
\]
and so $2^k < 3^{k-1}$. 
Note that between two consecutive integral powers of 2 there is at most one 
integral power of 3.  

\bt                    \label{NoGPunit:theorem:A3}
If $k \geq 3$, then 
\[
A_3^{(k)} = 2^{k-1}.
\]
\et

\begin{proof}
Let $G_1^{(k)} = (1/2^{k-1},1]$.
If 
\[
 \frac{1}{2^k} < x \leq \frac{1}{2^{k-1}} 
\]
then 
\[
\frac{1}{2^{k-1}} = \frac{2}{2^k} < 2x < 2^2 x < \cdots < 2^{k-1} x \leq \frac{2^{k-1} }{2^{k-1} } =1
\]
and so $\{2^ix:i=1,2,\ldots, k-1\}\subseteq G_1^{(k)}$, 
that is, $x$ is $k$-bad with respect to $G_1^{(k)}$.

If $x = 1/2^k$, then $2x = 1/2^{k-1} \notin G_1^{(k)}$.
If $r \geq 3$, then 
\[
r^{k-1} x \geq \frac{3^{k-1}}{2^k} > 1
\]
and so $r^{k-1} x  \notin G_1^{(k)}$.  
Therefore, $1/2^k$ is $k$-good with respect to $G_1^{(k)}$, 
and  $A_3^{(k)}  = 2^k$.
\end{proof}

\bt                    \label{NoGPunit:theorem:A4}
Let  $k \geq 3$.  
If there is no integral power of 3 between $2^{k-1} $ and $ 2^k$, then 
\[
A_4^{(k)}  = 
3^{k-1}.
\]
If there is an integral power of 3 between $2^{k-1} $ and $ 2^k$, 
and if $\ell$ is the positive integer such that  
\beq   \label{NoGPunit:2-3intertwine} 
2^{k-1} < 3^{\ell} < 2^k
\eeq
then $2 \leq \ell \leq k-2$  and 
\[
A_4^{(k)}  = 2^k 3^{k-1-\ell} = 3^{k-1}\left(\frac{2^k}{3^{\ell}}\right).
\]
\et

Inequality~\eqref{NoGPunit:2-3intertwine} is equivalent to 
$1 < 2^k/3^{\ell} < 2$.

For positive integers $k$, the following are equivalent: 
\benum
\item[(i)]
There is an integral power of 3 between $2^{k-1}$ and $2^k$.  
\item[(ii)]
The fractional part of $k\log_3 2$ is less than $\log_3 2$.
\item[(iii)]
$k$ is in the set $\{ [\ell \log_2 3]+1 : \ell = 1, 2, \ldots \}$. 
Thus, the formula for $A_4^{(k)}$ depends on diophantine properties of logarithms.
\eenum

\begin{proof}
We have $A_3^{(k)} = 1/2^k$ by Theorem~\ref{NoGPunit:theorem:A3}.  
Let
\[
 \frac{1}{4^{k-1}} < x \leq \frac{1}{2^k}.
\]
For $r\geq 4$ we have  
\[
x r^{k-1}  \geq x 4^{k-1} >  \frac{4^{k-1}}{4^{k-1}}  = 1
\]
and so $\{x r^i:i=1,2,\ldots, k-1 \} \not\subseteq  G_1^{(k)} $.  

With $r=2$ we have  
\[
x 2^{k-1} >  \frac{2^{k-1}}{4^{k-1}}  = \frac{1}{2^{k-1}}.
\]
Let $j$ be the smallest integer such that $x 2^j  >  1/2^{k-1}$.  
Then $j \leq k-1$.  Because $2 x \leq  2/2^k = 1/2^{k-1}$, it follows that $j \geq 2$.  
If 
\[
x 2^{j-1} \leq \frac{1}{2^k}
\]
then
\[
x 2^j  \leq \frac{1}{2^{k-1}} <  x 2^j 
\]
which is absurd.  Therefore,
\[
 \frac{1}{2^k} < x 2^{j-1} \leq \frac{1}{2^{k-1}} 
\]
and and so $\{ x 2^i:i=1,2,\ldots, k-1 \} \not\subseteq  G_1^{(k)} $.  

The remaining case is the ratio $r=3$ and the geometric progression 
$\{ x 3^i :i=1,2,\ldots, k-1 \}$.
If $x > 1/3^{k-1}$, then $x 3^{k-1} > 1$ and $x$ is good with respect to $G_1^{(k)} $.  
Therefore, $A_4^{(k)}  \geq 3^{k-1}$.  

Let $x = 1/3^{k-1}$.  
If there exists  $j \in \{1,2,\ldots, k-1\}$ such that 
\[
 \frac{1}{2^k}  <  x 3^j = \frac{3^j}{3^{k-1}}  < \frac{1}{2^{k-1}} 
 \]
then
\[
2^{k-1} < 3^{k-1-j} < 2^k.
\]
Thus, if there is no power of 3 between $2^{k-1}$ and $2^k$, then for all 
$i \in \{ 1,2,\ldots, k-1\}$, either $x < 3^ix < 1/2^k$ or $1/2^{k-1} < 3^i x \leq 1$.  
Thus, $1/3^{k-1}$ is $k$-bad, and $A_4^{(k)}  =3^{k-1}$.

Suppose that there is a power of 3 between $2^{k-1}$ and $2^k$, 
and that $\ell$ is the unique positive integer 
that satisfies~\eqref{NoGPunit:2-3intertwine}.
We observe that $k \geq 4$ because there is no power of 3 between $2^2=4$ 
and $2^3=8$, and that $2 \leq \ell \leq k-2$ 
because $2^3 < 3^2 \leq 3^{\ell}$ and $2^{k-1} < 3^{\ell} \leq 3^{k-2}$.  
Let 
\[
j = k-1 -\ell.
\]
Then $1 \leq j \leq k-3$.  
For $k \geq 4$ we have 
\[
\left(\frac{4}{3}\right)^{k-1} \geq \left(\frac{4}{3}\right)^3 > 2
\]
and so 
\[
\left(\frac{2}{3}\right)^{k-1} > \frac{1}{2^{k-2}}.
\]
Let  
\[
x_0 = \frac{1}{2^k3^j} = \frac{3^{\ell}}{2^k3^{k-1}} > \frac{2^{k-1}}{2^k3^{k-1}} 
=  \frac{1}{2^k} \left(\frac{2}{3}\right)^{k-1} 
>  \frac{1}{4^{k-1}}.
\]
If  
\[
x_0 < x \leq \frac{1}{3^{k-1}}
\]
then 
\[
\frac{1}{2^k} = x_0 3^j  < x 3^j  \leq  \frac{3^j}{3^{k-1}} 
= \frac{1}{3^{\ell} } < \frac{1}{2^{k-1}}
\]
and so  $x$ is good with respect to $G_1^{(k)} $.  

It remains only to prove that $x_0$ is bad.  
If $1 \leq i \leq j$, then 
\[
x_0 <  x_0 3^i \leq x_03^j  = \frac{1}{2^k}.
\]
If $j+1 \leq i \leq k-1$, then 
\[
\frac{1}{2^{k-1}} < \frac{3}{2^k} = 3^{j+1}x_0 \leq 3^ix_0 \leq 3^{k-1}x_0 
=  \frac{3^{\ell}}{2^k} < 1.
\]
Thus,  $x_0 = 1/(2^k3^j)$ is bad and $A_4^{(k)}  = 2^k3^j$.  
This completes the proof.  
\end{proof}

\section{Integer sequences with no $k$-term geometric progression}
If $a$ and $b$ are  real numbers with $a \leq b$, then the number of integers in the interval $(a,b]$ is $b-a + \theta$ with $|\theta| < 1$.

Recall that, for positive integers $k$ and $n$, the arithmetic function 
$g_k(n)$ denotes the cardinality of the largest subset of the set $\{1,2,3,\ldots, n\}$ 
that contains no integer geometric progression of length $k$ with integer ratio.

\bt       \label{NoGPunit:theorem:LowerBound}
Let $k \geq 3$, and let $\left(A_i^{(k)} \right)_{i=1}^{\infty}$ be the strictly increasing 
sequence of positive integers constructed in 
Theorem~\ref{NoGPunit:theorem:constructHarmonic}.
Then
\[
\gamma_k =  \liminf_{n\rightarrow \infty} \frac{ g_k(n)}{n} \geq 
  \sum_{i=1}^{\infty} \left(\frac{1}{A_{2i-1}^{(k)} } - \frac{1}{A_{2i}^{(k)} } \right).  
\]
In particular, 
\[
\gamma_k  \geq 1 - \frac{1}{2^k}  - \frac{1}{3^k}.
\] 
\et

\begin{proof}
For every positive integer $h$, the set 
\[
G_h^{(k)}  = \bigcup_{i=1}^h \left(\frac{1}{A_{2i}^{(k)} }, \frac{1}{A_{2i-1}^{(k)} }     \right]
\]
is a $k$-good subset of $(0,1]$.  
For every positive integer $n$, the dilated set 
\[
n \ast G_h^{(k)}  = n \ast \bigcup_{i=1}^h \left(\frac{1}{A_{2i}}, \frac{1}{A_{2i-1}}   \right] 
=  \bigcup_{i=1}^h \left(\frac{n}{A_{2i}}, \frac{n}{A_{2i-1}}   \right]
\]
is a disjoint union of intervals, and so
\begin{align*}
\left| (n \ast G_h^{(k)} ) \cap \N \right| 
& =  \sum_{i=1}^h  \left| \left(\frac{n}{A_{2i}}, \frac{n}{A_{2i-1}} \right] \cap \N \right| \\
& =  n \sum_{i=1}^h \left(\frac{1}{A_{2i-1}} - \frac{1}{A_{2i}} \right) + \theta_h 
\end{align*}
with $|\theta_h | < h$.
Because the dilation of a $k$-good set is $k$-good, 
and a subset of a $k$-good set is $k$-good, it follows that 
$(n \ast G_h^{(k)} ) \cap \N$ is a $k$-good set of positive integers.  
Moreover, $A_1 = 1$ implies that $(n \ast G_h^{(k)} ) \cap \N$ is a subset 
of $\{1,2,\ldots, n\}$.  Therefore, 
\[
\left| (n \ast G_h^{(k)} ) \cap \N \right| \leq g_k(n)
\]
and so 
\beq               \label{NoGPunit:gammaLowerBound}
 \sum_{i=1}^h \left(\frac{1}{A_{2i-1}} - \frac{1}{A_{2i}} \right) 
 = \lim_{n\rightarrow \infty} \frac{ \left| (n \ast G_h^{(k)} ) \cap \N \right| }{n} 
 \leq \liminf_{n\rightarrow \infty} \frac{ g_k(n)}{n}.
\eeq
This inequality holds for all $h \in \N$, and so 
\[
 \sum_{i=1}^{\infty} \left(\frac{1}{A_{2i-1}} - \frac{1}{A_{2i}} \right) 
 \leq \liminf_{n\rightarrow \infty} \frac{ g_k(n)}{n} = \gamma_k.
\]
Applying inequality~\eqref{NoGPunit:gammaLowerBound} with $h = 2$ and 
the values for $A_3^{(k)}$ and $A_4^{(k)}$ computed in
Theorems~\ref{NoGPunit:theorem:A3} and~\ref{NoGPunit:theorem:A4}, 
we obtain
\[
\gamma_k \geq \left(1 -  \frac{1}{2^{k-1}} \right) + \left( \frac{1}{2^k}  - \frac{1}{3^k} \right)
= 1 - \frac{1}{2^k}  - \frac{1}{3^k}.
\]
This completes the proof.
\end{proof}

It is a finite calculation to determine explicit values of the integers $A_i^{(k)}$ for small values of $i$ and $k$.  Table 1 contains all values of $A_i^{(k)}$ 
for $3 \leq k \leq 9$ with $A_i^{(k)} < 10^6$.
Applying inequality~\eqref{NoGPunit:gammaLowerBound} in Theorem~\ref{NoGPunit:theorem:LowerBound}, we can use 
these values to get lower bounds for $\gamma_k$ that improve results 
obtained previously by Rankin~\cite{rank60} and Riddell~\cite{ridd69}.  
For $k=3$, McNew~\cite{mcne13} has the current best lower bound.
Related results have been obtained by Brown and Gordon~\cite{brow-gord96},  
Beiglb{\"o}ck, Bergelson, Hindman, and Strauss~\cite{beig-berg-hind-stra06},
and Nathanson and O'Bryant~\cite{nath-obry13,nath-obry14}.

The following table records upper  and lower bounds for $\gamma_k$.

\begin{table}[h]\begin{center}  \footnotesize
\begin{tabular}{|c|cccc|c|ccc|} \hline
      & \multicolumn{4}{|c|}{Lower bounds on $\gamma_k$} & & \multicolumn{3}{|c|}{Upper bounds on $\gamma_k$} \\
  $k$ &  Rankin    & Riddell       & This paper & McNew      &        k               & McNew     & From $r_k$ &  Riddell     \\ \hline
  3   &  0.719 745 &                  & 0.815 870 & 0.818 410 & 3 & 0.819222 & 0.846 376  &  0.857 143 \\
  4   &  0.862 601 & 0.895 283 & 0.919 818 &                  & 4&                 & 0.928 874  &  0.933 334 \\
  5   &  0.931 652 & 0.958 056 & 0.963 737 &                  & 5 &                 & 0.967 742  &  0.967 742 \\
  6   &  0.966 324 & 0.980 371 & 0.982 877 &                  & 6 &                 & 0.983 871  &  0.984 126 \\
  7   &  0.983 438 & 0.991 159 & 0.991 805 &                  & 7 &                 &                   &  0.992 126 \\
  8   &  0.991 841 & 0.995 717 & 0.995 913 &                  & 8 &                 &                   &  0.996 079 \\
  9   &  0.995 969 & 0.997 939 & 0.998 012 &                  & 9 &                 &                   &  0.998 044 \\
  \hline
\end{tabular} 
\end{center}\end{table}

\def\cprime{$'$} \def\cprime{$'$} \def\cprime{$'$}
\providecommand{\bysame}{\leavevmode\hbox to3em{\hrulefill}\thinspace}
\providecommand{\MR}{\relax\ifhmode\unskip\space\fi MR }
\providecommand{\MRhref}[2]{%
  \href{http://www.ams.org/mathscinet-getitem?mr=#1}{#2}
}
\providecommand{\href}[2]{#2}

\begin{table}    \label{tableAik}
\begin{center} \footnotesize
\[
\begin{array}{c|ccccccc}
   & \multicolumn{7}{c}{k} \\
 i & 3 & 4 & 5 & 6 & 7 & 8 & 9 \\ \hline
 1 & 1 & 1 & 1 & 1 & 1 & 1 & 1 \\
 2 & 4 & 8 & 16 & 32 & 64 & 128 & 256 \\
 3 & 8 & 16 & 32 & 64 & 128 & 256 & 512 \\
 4 & 9 & 48 & 96 & 243 & 1152 & 2304 & 6561 \\
 5 & 12 & 200 & 144 & 288 & 1728 & 3456 & 6912 \\
 6 & 24 & 216 & 576 & 576 & 8192 & 16384 & 13824 \\
 7 & 27 & 288 & 4032 & 729 & 28800 & 32768 & 19683 \\
 8 & 32 & 1200 & 4096 & 1152 & 172800 & 163840 & 131072 \\
 9 & 36 & 1296 & 4608 & 2048 & 248832 & 288000 & 221184 \\
 10 & 40 & 1400 & 32256 & 3645 & 307328 & 331776 & 492075 \\
 11 & 45 & 1512 & 32768 & 4000 & 395136 & 497664 & 655360 \\
 12 & 48 & 1600 & 36288 & 10240 &           & 884736 &           \\
 13 & 2208 & 1728 & 36864 & 20736 &           & 995328 &           \\
 14 & 2209 & 1800 & 40320 & 21952 &           &           &           \\
 15 & 2256 & 1944 & 40960 & 92160 &           &           &           \\
 16 & 8832 & 2000 & 41472 & 100000 &           &           &           \\
 17 & 8836 & 62400 & 129600 & 102400 &           &           &           \\
 18 & 9024 & 63936 & 131072 & 207360 &           &           &           \\
 19 & 17664 & 73800 & 147456 & 219520 &           &           &           \\
 20 & 17672 & 74088 & 157216 & 518400 &           &           &           \\
 21 & 18048 & 75600 & 166464 & 548800 &           &           &           \\
 22 & 19872 & 79704 &           & 921600 &           &           &           \\
 23 & 19881 & 80688 &           &           &           &           &           \\
 24 & 20304 & 81648 &           &           &           &           &           \\
 25 & 26496 & 88000 &           &           &           &           &           \\
 26 & 26508 & 499200 &           &           &           &           &           \\
 27 & 27072 & 511488 &           &           &           &           &           \\
 28 & 52992 & 590400 &           &           &           &           &           \\
 29 & 53016 & 592704 &           &           &           &           &           \\
 30 & 54144 & 604800 &           &           &           &           &           \\
 31 & 59616 & 637632 &           &           &           &           &           \\
 32 & 59643 & 645504 &           &           &           &           &           \\
 33 & 60912 & 653184 &           &           &           &           &           \\
 34 & 70656 & 704000 &           &           &           &           &           \\
 35 & 70688 & 998400 &           &           &           &           &           \\
 36 & 72192 &           &           &           &           &           &           \\
 37 & 79488 &           &           &           &           &           &           \\
 38 & 79524 &           &           &           &           &           &           \\
 39 & 81216 &           &           &           &           &           &           \\
 40 & 88320 &           &           &           &           &           &           \\
 41 & 88360 &           &           &           &           &           &           \\
 42 & 90240 &           &           &           &           &           &           \\
 43 & 99360 &           &           &           &           &           &           \\
 44 & 99405 &           &           &           &           &           &           \\
 45 & 101520 &           &           &           &           &           &           \\
 46 & 103776 &           &           &           &           &           &           \\
 47 & 103823 &           &           &           &           &           &           \\
 48 & 105984 &           &           &           &           &           &           \\
 49 & 106032 &           &           &           &           &           &           \\
 50 & 108192 &           &           &           &           &           &           \\
 51 & 108241 &           &           &           &           &           &           \\
 52 & 108288 &           &           &           &           &           &           \\
\end{array}
\]
\end{center}
\caption{For $3\leq k \leq 9$, the table contains all integers $A^{(k)}_i$ 
satisfying Theorem~\ref{NoGPunit:theorem:constructHarmonic} 
that are less than $10^6$.
These numbers  are sequences A235054-60 in the OEIS. \label{ATable}}
\end{table}

\end{document}